\documentclass[a4paper,11pt, reqno]{amsart}
\usepackage[left= 3.4cm,right = 3.4cm, bottom=4.5cm, top=3.5cm, footskip = 1.5cm]{geometry}

\title[Auxiliary Monge-Amp\`ere Equations in Orbifold Setting]{Auxiliary Monge-Amp\`ere Equations in Orbifold Setting -- a Mean-Value Inequality}

\author{Johannes Scheffler}
\thanks{University of Bayreuth, Germany; e-mail: johannes.scheffler@uni-bayreuth.de}

\usepackage[utf8]{inputenc}
\usepackage[ngerman, english]{babel}
\usepackage[autostyle=true,german=quotes]{csquotes}
\usepackage[T1]{fontenc}
\usepackage{graphicx, subfig}

\usepackage{lmodern}
\usepackage{color}
\usepackage[maxbibnames=99,style=alphabetic]{biblatex}  
\usepackage{hyperref,mathtools}
\usepackage[skip=10pt plus2pt]{parskip}

\usepackage{xurl}
\hypersetup{breaklinks=true}

\usepackage{verbatim} 

\usepackage{amsfonts}
\usepackage{amsmath}
\usepackage{amsthm} 

\usepackage{amssymb}
\usepackage[new]{old-arrows}
\usepackage{faktor}
\usepackage{tikz-cd} 

\usepackage{mathrsfs}

\usepackage{relsize}

\usepackage{ dsfont }

\newcommand{\RR}{\mathbb{R}}
\newcommand{\CC}{\mathbb{C}}

\newcommand{\NN}{\mathbb{N}}

\newcommand{\XX}{\mathscr{X}}
\newcommand{\VV}{\mathscr{V}}

\newcommand{\dbar}{\overline{\partial}}
\newcommand{\ddbar}{\partial \overline{\partial}}

\newcommand{\Cont}{\mathscr{C}}

\newcommand{\ii}{\mathrm{i}}

\DeclareMathOperator{\id}{id}

\DeclareMathOperator{\Psh}{Psh}

\DeclareMathOperator{\tr}{tr}
\DeclareMathOperator{\Vol}{Vol}
\DeclareMathOperator{\Ent}{Ent}

\usepackage{enumitem}

\let\oldtocsection=\tocsection

\let\oldtocsubsection=\tocsubsection

\renewcommand{\tocsection}[2]{\hspace{0em}\oldtocsection{#1}{#2}}
\renewcommand{\tocsubsection}[2]{\hspace{1em}\oldtocsubsection{#1}{#2}}

\begin{document}

	\theoremstyle{definition}
	\newtheorem{Definition}{Definition}
	\newtheorem{Beispiel}[Definition]{Example}
	\newtheorem{Construction}[Definition]{Construction}
	
	\theoremstyle{remark}
	\newtheorem{Bemerkung}[Definition]{Remark}
	
	\theoremstyle{plain}
	\newtheorem{Satz}[Definition]{Satz}
	\newtheorem{Folgerung}[Definition]{Corollary}
	\newtheorem{Lemma}[Definition]{Lemma}
	\newtheorem{Proposition}[Definition]{Proposition}
	\newtheorem{Theorem}[Definition]{Theorem}
	\newtheorem{Conjecture}[Definition]{Conjecture}
	\newtheorem*{THeorem}{Theorem} 
	\newtheorem*{SAtz}{Satz}

	\maketitle

	\pagenumbering{arabic} 
	
	\begin{abstract}
		In this note, we generalize the mean-value-inequality of Guo-Phong-Sturm (cf.\ Lemma 2 in \cite{GPS22}) to the setting of a compact Kähler orbifold. This shows that their reasoning is insensitive to quotient singularities. As we aim for a self-contained exposition, we generalize some fundamental results: Hörmander's and Tian's $\alpha$-invariant estimate (see Section \ref{sec_alpha_est}), Berman's approximation of the psh-envelope of a $(1,1)$-form in a Kähler class (see Section \ref{sec_approx_envelope}) and an $L^\infty$ estimate for this envelope by Guo-Phong-Tong-Wang (see Section \ref{sec_Linfty_est}).
	\end{abstract}
	
	\tableofcontents
	\enlargethispage{1\baselineskip}

	\section{Setup \& results}
	
	During the last years, the PDE approach to Monge-Ampère equations on compact Kähler manifolds developed by Guo-Phong-Tong led to important breakthroughs and new proofs of known results. For instance, estimates for the Green kernel in \cite{GPS22} and subsequently diameter estimates in \cite{GPSS24} have been obtained without assuming bounds on the Ricci curvature. An essential basis for these results is a mean-value inequality (see Lemma 2 in \cite{GPS22} and Theorem \ref{meanvalueineq}). This inequality has also been applied by Cao-Graf-Naumann-Peternell-P\u{a}un-Wu in the realm of Hermite-Einstein metrics on singular spaces (cf. \cite{CGNPPW23}). 
	
	Since non-smooth spaces are of high interest, e.g.\ in  the MMP, the question whether the method of Guo-Phong-Tong using auxiliary Monge-Amp\`ere equations generalizes to singular settings is very natural. By developing the theory needed for the mean-value inequality, we answer this question affirmatively in the case of Kähler orbifolds. 
	
	We will start by introducing the set-up. Afterwards, the main results are stated. Whenever possible, we will use the same notations as in the original proof in \cite{GPS22}. Thereby, we hope to improve the readability for those who are already familiar with the arguments in the smooth setting. Let $(\XX , \omega_\XX)$ be a compact Kähler orbifold of complex dimension $n$. We fix a locally finite orbifold atlas $\{(U_\alpha, G_\alpha, \pi_\alpha)\}$. Without loss of generality, the groups $G_\alpha$ act on $U_\alpha \subset \CC^n$ as regular linear maps. We understand $\omega_\XX$ to be a collection of $G_\alpha$-invariant closed, positive, real $(1,1)$-forms $\omega_{\XX,\alpha}$ on $U_\alpha$ satisfying the usual compatibility conditions for orbifolds.\par
	
	We suppose that there is a real, $d$-closed $(1,1)$-form $\chi$ on $\XX$. It defines a cohomology class $[\chi] \in H^{1,1}(\XX, \RR)$. We will use the following notions:
	
	\begin{Definition}
		A class $[\chi] \in H^{1,1}(\XX, \RR)$ is said to be \emph{nef} if for every $\varepsilon > 0$ we find a representative $\Tilde{\chi} \in [\chi]$ such that
		\begin{align}\label{eq_def_nef}
			\Tilde{\chi} > - \varepsilon \omega_\XX.
		\end{align}
		The class $[\chi]$ is said to be \emph{big} if 
		\begin{align*}
			\int_\XX \chi^n > 0.
		\end{align*}
	\end{Definition}
	
	\begin{Bemerkung}
		\begin{enumerate}
			\item The equation (\ref{eq_def_nef}) is understood to hold for the corresponding forms $\chi_\alpha$ and $\omega_{\XX,\alpha}$ on each $U_\alpha$. 
			\item Since the Theorem of Stokes holds in the orbifold setting, the definition of bigness is independent of the choice of representative.
		\end{enumerate}
	\end{Bemerkung}
	
	Assume from now on that $[\chi]$ is nef and big. Despite $\chi$ and $\hat{\omega}_t \coloneqq \chi + t \omega_\XX$ need not to be semi-positive, this implies as in the smooth setting that for every $t >0$ there is a Kähler form $\omega_t$ in the class $[\chi+t\omega_\XX]= [\hat{\omega}_t]$. Simply choose a representative $\Tilde{\chi} \in [\chi]$ such that $\Tilde{\chi} > - \frac{t}{2} \omega_\XX$.
	
	Now, we assume that the volume of $\XX$ with respect to $\omega_\XX$ is normalized as $V \coloneqq \int \omega_\XX^n = 1$, and denote for each $t$ the volume with respect to $\omega_t$ by $V_t \coloneqq \int \omega_t^n$. By the big- and nefness assumption on $\chi$, $V_t$ is bounded from below by $V_0 \coloneqq \int \chi^n >0$. We need the function
	\begin{align*}
		F_t \coloneqq \log \left ( \frac{1}{V_t}\frac{ \omega_t^n}{ \omega_\XX^n}\right)
	\end{align*}
	to compare the involved volume forms. Clearly, this definition works in charts $U_\alpha$ and descends to $\XX$ as all objects involved are $G_\alpha$ invariant.
	
	We fix $p > n$ and denote the $p$-th entropy of $\omega_t$ as in the smooth case by
	\begin{align*}
		\Ent_p(\omega_t) \coloneqq \frac{1}{V_t} \int_\XX |F_t|^p \omega_t^n =  \int_\XX |F_t|^p e^{F_t} \omega_\XX^n.
	\end{align*}
	We have to assume that the entropy of the family $\{\omega_t\}$ is uniformly bounded by a constant $N> 0$. This assumption turns out to be fulfilled in applications (cf.\ e.g.\ \cite{CGNPPW23}) and is clearly less restrictive than assumptions on the Ricci curvature involving two derivatives more. We would like to mention at this point that replacing assumptions on derivatives by integral bounds goes back to the work of S. Kołodziej (cf.\ \cite{Kolodziej98_CplxMA}).
	
	By the $\partial \dbar$-Lemma for Kähler orbifolds, for every $t$, there is $\varphi_t$ such that
	\begin{align*}
		\omega_t = \chi + t\omega_\XX + i \partial \dbar \varphi_t.
	\end{align*}
	The function $\varphi_t$ becomes unique if we demand $\sup_\XX \varphi_t = 0$, and solves the Monge-Ampère equation
	\begin{align*}
		(\chi + t\omega_\XX + i \partial \dbar \varphi_t)^n = V_t e^{F_t} \omega_\XX^n.
	\end{align*}
	
	To compensate $\chi$'s and $\hat{\omega}_t$'s possible negativity, we need to work with the \emph{plurisubharmonic envelope} associated to $\hat{\omega}_t$ defined by
	\begin{align*}
		\VV_t \coloneqq \sup \{u \colon \XX \longrightarrow \RR \text{ usc } \mid u \leq 0, \quad \hat{\omega}_t + \ii \partial \dbar u \geq 0 \}.
	\end{align*}
	
	We start now with our first main result, generalizing the uniform $L^\infty$ estimate for the family $\{\varphi_t\}$ (cf.\ \cite[Thm 1 (b)]{GPTW24}) to the orbifold setting.
	
	\begin{Theorem} \label{Linfty_estimate}
		There is a uniform constant $C > 0$ depending only on $n$, $p$, $(\XX, \omega_\XX)$, $\chi$, and $N$ such that for all $t \in (0,1]$ we have
		\begin{align*}
			\sup_\XX |\varphi_t - \VV_t | \leq C.
		\end{align*}
	\end{Theorem}
	
	The proof in Section \ref{sec_Linfty_est} follows the same structure as in the smooth setting. However, we provide some additional explanations and a couple of simplifications. In particular, we generalize the ingredients of the proof and streamline them for our purposes. The first of these is an approximation result for the psh-envelope by R.\ Berman (cf. \cite{Berman2019}) depending on J.-P.\ Demailly's approximation by functions with analytic singularities (cf.\ \cite{Dem1992_RegCurIntTheo} \cite{Wu23}).
	
	\begin{Proposition} \label{approx_psh_envelope}
		For any $\beta >  0$ the Monge-Ampère equation
		\begin{align*}
			(\hat{\omega}_t + \ii \partial \dbar u_\beta )^n = V_t e^{\beta u_\beta} \omega_\XX^n
		\end{align*}
		admits an unique, orbifold-smooth solution $u_\beta$ which converges uniformly to the psh-envelope $\VV_t$.
	\end{Proposition}
	
	Moreover, as second ingredient, we need Hörmander's and Tian's $\alpha$-invariant estimate for $\omega_\XX$-psh functions. Section \ref{sec_alpha_est} is devoted to its proof for compact Kähler orbifolds. Then, Theorem \ref{Linfty_estimate} and the $\alpha$-invariant estimate are used in Section \ref{sec_meanvalue_ineq} to prove the second main result, the aforementioned mean-value-inequality which was established by \cite{GPS22} in the smooth setting:
	
	\begin{Theorem} \label{meanvalueineq}
		Let $a > 0$ be fixed. There is a constant $C > 0$ depending only on $n$, $p$, $\chi$, $\omega_\XX$, $N$, and $a$ such that all $v \in L^1(\XX, \omega_t)$ which 
		\begin{enumerate}
			\item are twice continuously differentiable on $\Omega_{-c} \coloneqq \{v > -c\}$ for some $c > 0$,
			\item satisfy $\int_\XX v \omega_t^n = 0$, and
			\item whose Laplacian with respect to $\omega_t$ is bounded in $\Omega_0$ from below by
			\begin{align}
				\Delta_{\omega_t} v \geq -a
			\end{align}
		\end{enumerate} 
		are bounded from above by
		\begin{align}
			\sup_\XX v \leq C(1 + \|v\|_{L^1(\XX, \omega_t)}).
		\end{align}
	\end{Theorem}
	
	\textbf{Acknowledgements.} I would like to thank my advisor Mihai P\u{a}un for pointing out the relevance of the mean-value inequality, proposing its generalization to orbifolds, and guiding me throughout this project. Also, I am grateful to Philipp Naumann for discussions about orbifolds and to Duong H.\ Phong and Freid Tong for explanations concerning their method of proof. Finally, I would like to acknowledge the financial support received by the Studienstiftung des deutschen Volkes and the Marianne-Plehn Programm.
	
	\bigskip
	
	
	
	
	\section{Hörmander's and Tian's \texorpdfstring{$\alpha$}{alpha}-invariant estimate on compact Kähler orbifolds} \label{sec_alpha_est}
	
	Throughout this section, let $(\XX , \omega_\XX)$ be a compact orbifold endowed with a fixed orbifold Kähler metric. We want to generalize Hörmander's and Tian's estimate for $\omega_\XX$-psh functions, more precisely functions in the space
	\begin{align*}
		P(\XX, \omega_\XX) \coloneqq \left \{ \varphi \in \Cont^2(\XX, \RR) \mid \omega_\XX + \ii \partial \dbar \varphi \geq 0 , \quad \sup_\XX \varphi = 0 \right \} \subset \Psh(\XX, \omega_\XX).
	\end{align*}
	We have the same result as stated by G.\ Tian in Proposition 2.1 of \cite{Tian87_KEMetrics} in the smooth setting:
	
	\begin{Proposition} \label{alpha_inv_est}
		There exist positive constants $\alpha > 0$ and $C > 0$ which only depend on $(\XX, \omega_\XX)$ such that for all $\varphi \in P(\XX, \omega_\XX)$ we have
		\begin{align*}
			\int_\XX e^{-\alpha \varphi} \omega_\XX^n \leq C.
		\end{align*}
	\end{Proposition}
	
	We will use the same approach for the proof as in \cite{Tian87_KEMetrics}, in particular the same local result which is essentially due to L.\ Hörmander (cf.\ Theorem 4.4.5 in \cite{Hor73_SCV}):
	
	\begin{Lemma} \label{local_hormander_lemma}
		Let $B_R(0)$ denote the ball of radius $R$ around $0$ in $\CC^n$, and let $\lambda > 0$ be fixed. Then there is a constant $C > 0$ depending only on $n, \lambda$, and $R$ such that for any plurisubharmonic function $\psi \in \Psh(B_R(0))$ satisfying $\psi \leq 0$ in $B_R(0)$ and $\psi(0) \geq -1$ we have
		\begin{align*}
			\int_{|z| < r} e^{-\lambda \psi(z)} dz \leq C
		\end{align*}
		for any $r < R e^{-\lambda/2}$.
	\end{Lemma}
	
	\begin{proof}[Proof of Proposition \ref{alpha_inv_est}.] 
		As in the original proof, we want to estimate the supremum of $\varphi \in P(\XX, \omega_\XX)$ from below on small balls, of course independently of the chosen $\varphi$. For this, consider the Green function $G$ of the Laplacian $\Delta$ with respect to the metric $\omega_\XX$ such that its infimum is $\inf_{\XX \times \XX} G (x,y) = 0$ (cf.\ \cite{Faulk2019_thesis, Chiang1990_HarmonicMapsOrbifolds}). We have the usual Green-Riesz representation formula
		\begin{align*}
			\varphi(x) = \frac{1}{V} \int_\XX \varphi (y) \omega_\XX^n (y) - \int_\XX G(x,y) \Delta \varphi (y) \omega_\XX^n (y).
		\end{align*}
		By assumption, $\omega_\XX + \ii \partial \dbar \varphi \geq 0$. Thus, taking the trace gives 
		\begin{align*}
			\Delta \varphi = \tr_{\omega_\XX} \ii \partial \dbar \varphi \geq - \tr_{\omega_\XX} \omega_\XX = -n
		\end{align*}
		so that 
		\begin{align*}
			\varphi(x) \leq  \frac{1}{V} \int_\XX \varphi (y) \omega_\XX^n (y) + n \int_\XX G(x,y) \omega_\XX^n (y).
		\end{align*}
		Taking the supremum on both sides yields
		\begin{align*}
			0 = \sup_\XX \varphi(x) \leq  \frac{1}{V} \int_\XX \varphi (y) \omega_\XX^n (y) + n \sup_\XX \int_\XX G(x,y) \omega_\XX^n (y). 
		\end{align*}
		Denote the supremum of the $L^1$-norm of $G$ by $C_1 > 0$, depending only on $(\XX, \omega_\XX)$. Then the average of $\varphi$ can be estimated from below by
		\begin{align} \label{Est_avg_psh_fct}
			\frac{1}{V} \int_\XX \varphi (y) \omega_\XX^n (y) \geq -n C_1.
		\end{align}
		Now, we fix a cover of $\XX$ as follows: around each point $x \in \XX$, choose an orbifold chart $(U_{\alpha, x}, G_{\alpha, x}, \pi_{\alpha, x})$ and a geodesic ball $B(p_x, r_x) \subset U_{\alpha, x}$ of radius $r_x >0$ around a preimage $p \in \pi^{-1}_{\alpha, x}(x)$ of $x$. Since $\XX$ is compact, there are finitely many $x_1, \dots , x_N$ such that the images of the corresponding balls with radius $r_{x_i}/4$ cover $\XX$:
		\begin{align*}
			\XX = \bigcup_{i=1}^N \pi_{x_i}(B(p_{x_i}, r_{x_i}/4)).
		\end{align*}
		From now on, we index with $i$ instead of $x_i$. 
		Our goal is to apply Lemma \ref{local_hormander_lemma} in each of these balls. Using that $\varphi \leq 0$ and (\ref{Est_avg_psh_fct}), we can estimate the integral of $\varphi$ on the image of each ball by
		\begin{align*}
			\int_{\pi_i(B(p_i, r_i/4))} \varphi (y) \omega_\XX^n (y) \geq \int_\XX \varphi (y) \omega_\XX^n (y) \geq -n C_1 V.
		\end{align*}
		Because the supremum of a function must be at least its average, this implies the existence of a point $q_i \in B(p_i, r_i/4)$ in each ball such that 
		\begin{align*}
			\varphi( \pi_i (q_i) ) \geq \frac{-nVC_1}{\Vol (\pi_i(B(p_i, r_i/4))) }.
		\end{align*}
		On each ball, let $\psi_i$ be the potential of the local representative $\omega_{\XX,i}$ of the Kähler form with the normalization that $\psi_i(q_i) = 0$. We denote the supremum of all potentials $\psi_i$ by 
		\begin{align*}
			C_2 \coloneqq \max_i \sup_{B(p_i, 7r_i/8)} |\psi_i |.
		\end{align*}
		We fix $\lambda \coloneqq \ln 5/4$ such that $r_i/2 = \frac{4}{8}r_i < \frac{5}{8} r_i e^{-\lambda/2}$ and we can apply Lemma \ref{local_hormander_lemma} with these radii. Then we set $V_{\min} \coloneqq \min \Vol (\pi_i(B(p_i, r_i/4)))$ and 
		\begin{align*}
			\alpha \coloneqq \frac{\lambda}{C_2 + nVC_1/V_{\min} }
		\end{align*}
		which depends only on properties of $(\XX, \omega_\XX)$ but not on $\varphi$. The function $\alpha/ \lambda \, (\psi_i (x) + \varphi(x) - C_2) = \frac{1}{C_2 + nVC_1/V_{\min} } (\psi_i (x) + \varphi(x) - C_2)$ satisfies the requirements of Lemma \ref{local_hormander_lemma} on $B(q_i , 5r_i/8)$:
		\begin{enumerate}
			\item it is psh because $\ii \partial \dbar (\psi_i  + \varphi) = \omega_{\XX, i} + \ii \partial \dbar \varphi \geq 0 $ since $\varphi$ is $\omega_\XX$-psh,
			\item it is non-positive because $\varphi \leq 0$ and $\psi_i \leq C_2$ on the larger ball $B(p_i, 7r_i/8)$, and
			\item it is at least $-1$ at the point $q_i$ since $\psi_i (q_i)=0$ and $\varphi (q_i) \geq -nV C_1/V_{\min}$ by the choice of $q_i$.
		\end{enumerate}
		Thus, the local Lemma yields 
		\begin{align*}
			\int_{B(q_i, r_i/2)} e^{-\alpha  (\psi_i  + \varphi - C_2)} \omega_{\XX,i}^n \leq C_i.
		\end{align*}
		We remark that the balls of radius $r_i/4$ around $p_i$ are contained in those of radius $r_i/2$ around $q_i$, and that $-\alpha  (\psi_i - C_2) \geq 0$ so that
		\begin{align*}
			\int_{\pi_i(B(p_i, r_i/4))} e^{-\alpha \varphi} \omega_{\XX}^n \leq \int_{B(p_i, r_i/4)} e^{-\alpha \varphi} \omega_{\XX,i}^n \leq \int_{B(q_i, r_i/2)} e^{-\alpha  (\psi_i  + \varphi - C_2)} \omega_{\XX,i}^n \leq C_i
		\end{align*}
		where in the first inequality the definition of local integration on orbifolds is used.
		Because the balls around $p_i$ of radius $r_i/4$ cover $\XX$, we conclude that
		\begin{align*}
			\int_\XX e^{-\alpha \varphi} \omega_{\XX}^n \leq \sum_{i=1}^N C_i \eqqcolon C.
		\end{align*}
	\end{proof}

	
	
	
	\section{Approximation of the psh-envelope in the orbifold setting} \label{sec_approx_envelope}
	In this section we want to explain Proposition \ref{approx_psh_envelope}. However, we will work in a slightly more general setting (similar to the Kähler case in \cite{Berman2019}, cf.\ section 2.1 there). Let $\XX$ be a compact orbifold, endowed with a fixed volume form $dV$. We assume there is an orbifold smooth but not necessarily positive representative $\theta$ of a Kähler class $[\theta] \in H^{1,1}(\XX, \RR)$. 
	
	\begin{Lemma} 
		The Monge-Ampère equation
		\begin{align} \label{psh_env_approx_MAeq}
			(\theta + \ii \partial \dbar u_\beta )^n = ce^{\beta u_\beta} dV 
		\end{align}
		admits an unique, orbifold smooth solution $u_\beta$. Here, $c$ is a normalizing constant equal to the ratio of $\int_\XX \theta^n$ and $\int_\XX dV$.
	\end{Lemma}
	
	\begin{proof}
		We reduce the problem to a Monge-Ampère equation whose solvability is well-known. By the $\ddbar$-Lemma there is a smooth function $\varphi$ such that $\theta + \ii \ddbar \varphi \eqqcolon \omega$ is a Kähler form.  We set $F\coloneqq \log (dV/\omega^n) + \log c' + \beta\varphi$ where $c'$ is a normalizing constant and solve
		\begin{align}
			((\beta\omega) + \ii \ddbar u )^n = e^{F+u} (\beta\omega)^n.
		\end{align}
		This Monge-Ampère equation admits an unique, orbifold-smooth solution (see e.g. Section 2.5 in \cite{Faulk2019_thesis} or \cite{Faulk2019_YauOrbifolds}). We use the solution to define
		\begin{align*}
			u_\beta \coloneqq \frac{u}{\beta} + \varphi.
		\end{align*}
		Then we have
		\begin{align*}
			(\theta + \ii \ddbar u_\beta)^n &= \frac{1}{\beta^n} ((\beta\omega) + \ii \ddbar u)^n\\ 
			&= \frac{1}{\beta^n} e^{F+u} (\beta\omega)^n \\
			&= \frac{dV}{\omega^n} c e^{\beta u_\beta} \omega^n
		\end{align*}
		establishing the Lemma.
	\end{proof}
	
	We denote by
	\begin{align} \label{def_psh_envelope}
		\VV_\theta \coloneqq \sup \{ u \leq 0 \mid u \text{ usc.\ and } \theta + \ii \ddbar u \geq 0 \}
	\end{align}
	the plurisubharmonic envelope of $\theta$. The supremum is taken over all non-positive, $\theta$-psh functions. However, by approximation with analytic singularities à la Demailly, it is enough to consider only smooth, $\theta$-psh functions:
	
	\begin{Lemma} \label{psh_env_is_smooth_psh_env}
		For a smooth representative $\theta$ of a Kähler class, the psh-envelope coincides with 
		\begin{align*}
			\VV_{\theta,\mathrm{sm}} \coloneqq \sup \{ u \leq 0 \mid u \text{ smooth and } \theta + \ii \ddbar u \geq 0 \}.
		\end{align*}
	\end{Lemma}
	
	\begin{proof}
		First, because the class $[\theta]$ is Kähler, by the $\ii \ddbar$-Lemma, there is an orbifold-smooth function $\varphi\leq 0$ such that $\omega \coloneqq \theta + \ii \ddbar \varphi > 0$ is a Kähler form. Therefore, $\VV_\theta$ is bounded from below by $\varphi$.\\
		Next, note that $\VV_\theta$ itself is $\theta$-psh. In general, when taking the supremum of psh functions, one has to take the so-called usc-regularization to obtain a psh-function. In our case however, the usc-regularized supremum is already an element of the set in (\ref{def_psh_envelope}).\\
		We apply now J.-P.\ Demailly's approximation by functions with analytic singularities (see \cite{Dem1992_RegCurIntTheo} Proposition 3.7 for the smooth and Theorem 6 in \cite{Wu23} for the orbifold case) to $\VV_\theta$ to obtain a sequence $(u_m)$ of functions with the following properties: $u_m \geq \VV_\theta$ and $u_m$ converges to $\VV_\theta$ pointwise and in $L^1$; the loss of positivity is bounded by $\ii \ddbar u_m \geq -\theta - \varepsilon_m \omega$ for a sequence $(\varepsilon_m)$ converging to zero. Note that quasi-psh functions with analytic singularities are smooth if they are bounded from below. Here, we find that $\VV_\theta$ and hence the $u_m$ are bounded from below. The approximations $u_m$ are consequently smooth!\\
		We define
		\begin{align*}
			u'_m \coloneqq (1-\varepsilon_m) u_m + \varepsilon_m \varphi - \sup_X \left( (1-\varepsilon_m) u_m + \varepsilon \varphi \right)
		\end{align*}
		to make the approximating sequence $\theta$-psh and non-positive. Indeed, we have
		\begin{align*}
			\theta + \ii \ddbar u'_m &= (1-\varepsilon_m) (\theta + \ii \ddbar u_m) + \varepsilon_m (\theta+  \ii \ddbar \varphi)\\ & \geq (1-\varepsilon_m)( - \varepsilon_m \omega) + \varepsilon_m \omega \\
			& \geq 0.
		\end{align*}
		By definition of $\VV_{\theta,\mathrm{sm}}$ we have for every $m$
		\begin{align*}
			\VV_{\theta,\mathrm{sm}} \geq u'_m \longrightarrow \VV_\theta
		\end{align*}
		for $m \rightarrow \infty$ which concludes the proof.
	\end{proof}
	
	We can now state the main result of this section.
	
	\begin{Proposition}
		There is a constant $C > 0$ such that 
		\begin{align*}
			|u_\beta - \VV_\theta | \leq \frac{C\log \beta}{\beta},
		\end{align*}
		in particular, the smooth functions $u_\beta$ approximate the psh-envelope uniformly.
	\end{Proposition}
	
	Note that in general, the constant $C$ may depend on the form $\theta$ and also on the choice of a Kähler representative in its class. However, for our application, the rate of convergence is not relevant.
	
	\begin{proof}
		It suffices to show the statement for $ \VV_{\theta,\mathrm{sm}}$ by Lemma \ref{psh_env_is_smooth_psh_env}. First, we establish the (simpler) bound of $u_\beta - \VV_{\theta,\mathrm{sm}}$ from above. Let $x_0 $ be the point where $u_\beta$ attains its maximum. Choose an orbifold chart $(U, G, \pi)$ around $x_0$ and a preimage $y_0 \in U$ of $x_0$. Then the smooth representative of $u_\beta$ on $U$ attains its maximum at $y_0$ and we have -- in terms of the representatives on the locally uniformizing system -- $\ii \ddbar u_\beta|_{y_0} \leq 0$ by the maximum principle. By the defining Monge-Ampère equation (\ref{psh_env_approx_MAeq}) it follows that  at $y_0$
		\begin{align*}
			\theta^n \geq e^{\beta u_\beta} dV
		\end{align*}
		which is equivalent to 
		\begin{align*}
			e^{\beta u_\beta} \leq \frac{\theta^n}{dV}.
		\end{align*}
		Because $u_\beta$ is maximal at $x_0$, on all of $\XX$ we have
		\begin{align*}
			u_\beta \leq \frac{1}{\beta} \log \sup_\XX \frac{\theta^n}{dV} \eqqcolon \frac{C_1}{\beta}.
		\end{align*}
		This implies that $u_\beta - C_1/\beta $ is smooth, $\theta$-psh, and non-positive and therefore bounded by $\VV_{\theta,\mathrm{sm}}$ from above.
		
		On the other hand, let $C_2 > 0$ be a constant such that $\omega^n \geq C_2 dV$ and let 
		\begin{align*}
			C_3 \coloneqq n + \max_{\beta \geq 2} \left\{0, \frac{\log C_2}{- \log \beta} \right\} > 0.
		\end{align*}
		For an arbitrary, smooth, non-positive, and $\theta$-psh function $v \in \Cont^\infty \cap \Psh(\XX, \theta)$ we define 
		\begin{align*}
			v_\beta \coloneqq \left( 1- \frac{1}{\beta}\right) v + \frac{1}{\beta} \varphi  - \frac{C_3 \log \beta}{ \beta}.
		\end{align*}
		Then we find that
		\begin{align*}
			(\theta + \ii \ddbar v_\beta )^n &= (\theta + (1-1/\beta) \ii \ddbar v + (1/\beta) \ii\ddbar \varphi )^n \geq 1/\beta^n (\theta + \ii\ddbar \varphi)^n 
		\end{align*}
		since $v$ is $\theta$-psh. By definition, $(\theta + \ii \ddbar \varphi )^n = \omega^n $ which can be estimated from below by $C_2 dV$. By the choice of $C_3$ we have that
		\begin{align*}
			\exp(- C_3 \log \beta) \leq \exp(- n \log \beta + \log C_2) = \frac{C_2}{\beta^n}.
		\end{align*}
		Combining the non-positivity of $v$ and $\varphi$ with these estimates yields
		\begin{align*}
			(\theta + \ii \ddbar v_\beta )^n &\geq \frac{C_2}{\beta^n} e^{(\beta-1)v + \varphi} dV \geq e^{\beta v_\beta } dV.
		\end{align*}
		This Monge-Ampère inequality implies that $v_\beta \leq u_\beta$ as follows. Since both, $u_\beta$ and $v_\beta$ are smooth, we can apply the maximum principle as before to $v_\beta - u_\beta$ and get $\ii \ddbar v_\beta \leq \ii \ddbar u_\beta$ at some point $x_0$. Thus,
		\begin{align*}
			e^{\beta v_\beta(x_0)} dV \leq (\theta + \ii \ddbar v_\beta )^n|_{x_0} \leq (\theta + \ii \ddbar u_\beta )^n|_{x_0} = e^{\beta u_\beta(x_0)} dV
		\end{align*}
		holds and therefore $\max_\XX (v_\beta - u_\beta) \leq 0$.\\
		In other words, 
		\begin{align*}
			\left( 1- \frac{1}{\beta}\right)v + \frac{1}{\beta} \varphi - \frac{C_3 \log \beta}{ \beta} \leq u_\beta
		\end{align*}
		holds for arbitrary smooth, non-positive, $\theta$-psh functions which we can thus replace by their pointwise supremum, $\VV_{\theta, \mathrm{sm}}$. This yields 
		\begin{align*}
			\VV_{\theta, \mathrm{sm}} - u_\beta &\leq \frac{1}{\beta}\VV_{\theta, \mathrm{sm}}  - \frac{1}{\beta} \varphi + \frac{C_3 \log \beta}{ \beta}\\
			&\leq \frac{C_4 \log \beta}{ \beta}
		\end{align*}
		because $\VV_{\theta, \mathrm{sm}} \leq 0$ and $\varphi$ is bounded. The proof is complete.
	\end{proof}

	
	
	
	\section{Mean-value inequality} \label{sec_meanvalue_ineq}
	This section is devoted to the proof of Theorem \ref{meanvalueineq} which goes along the same lines as in \cite{GPS22} and is very similar as the proof of Theorem \ref{Linfty_estimate} in the next section. Throughout this proof, we will assume Theorem \ref{Linfty_estimate} holds.
	
	First, we fix $0 < t \leq 1$ and $v\in L^1(\XX, \omega_t)$. We may assume that $v$ is not constantly zero and satisfies $\|v\|_{L^1} \leq V_0 = [\chi]^n  $. If not, we replace it by $V_0 v/\|v\|_{L^1}$ which is still $\Cont^2$ wherever non-negative, has mean-value zero on $\XX$ with respect to $\omega_t$, and its Laplacian is bounded from below by $-a$. Recall that we normalized the fixed reference metric $\omega_\XX$ to have volume $V=1$ for simplicity.
	
	We need to show that $\sup_\XX v \leq C$ where $C$ is indepent of $v$ and $t$. This is equivalent to proving that the function
	\begin{align*}
		\phi \colon \RR_{\geq 0} \longrightarrow \RR_{\geq 0}, \quad \quad s \longmapsto \frac{1}{V_t} \int_{\Omega_s} \omega_t^n = \int_{\Omega_s} e^{F_t} \omega_\XX^n
	\end{align*}
	vanishes for $s \geq C$. Here and in the following, $\Omega_s = \{v > s\}$ denotes the upper-level set of $v$. The main part of the proof is to establish for every $s\geq 0$ and every $r \geq 0$ the functional inequality 
	\begin{align} \label{fctl_ineq_phi}
		r \phi(s+r) \leq C' \phi (s)^{1+\delta_0}
	\end{align}
	where $\delta_0 \coloneqq (p-n)/(np) > 0$ and $C' >0$ depends on the same quantities as $C$. We prove first that a DeGiorgi-type iteration yields Theorem \ref{meanvalueineq} if (\ref{fctl_ineq_phi}) holds.
	
	\subsection{DeGiorgi-type iteration} \label{subsec_DeGiorgi}
	We first note that $\phi$ is not-increasing and decays at least at rate of $1/s$. On $\Omega_s$ for $s>0$ we have $v/s >1$ so that
	\begin{align*}
		\phi(s) < \int_{\Omega_s} \frac{v}{s} e^{F_t} \omega_X^n \leq \frac{1}{s} \int_{\Omega_0} v e^{F_t} \omega_X^n \leq  \frac{\|v\|_{L^1(\XX , \omega_t)}}{V_t s} < \frac{1}{s}
	\end{align*}
	by the assumption $\|v\|_{L^1(\XX , \omega_t)} \leq V_0 < V_t$. We choose $s_0 \coloneqq (2C')^{1/\delta_0}$ such that 
	\begin{align*}
		\phi (s_0)^{\delta_0} < \frac{1}{2C'}.
	\end{align*}
	Moreover, for $i \in \{0,1,2,\dots \}$ we define $r_i \coloneqq 2^{-i\delta_0}$ and we iterate (\ref{fctl_ineq_phi}) starting from $s_0$ with the $r_i$'s. We set
	\begin{align*}
		C \coloneqq s_0 + \sum_{i = 0}^{\infty} r_i = (2C')^{1/\delta_0} + \frac{1}{1-2^{-\delta_0}}
	\end{align*}
	which only depends on the quantities stated in Theorem \ref{meanvalueineq}.
	
	\textbf{Claim.} For every $m \geq 0$ we have
	\begin{align*}
		\phi\left(s_0 + \sum_{i=0}^m 2^{-i\delta_0} \right) \leq C'^{-1/\delta_0} 2^{-m-1-1/\delta_0}.
	\end{align*}
	
	If the claim holds, $\phi(s) = 0$ if $s \geq C$ so Theorem \ref{meanvalueineq} is established.
	
	\begin{proof}
		The proof proceeds by induction on $m$. If $m=0$ we use (\ref{fctl_ineq_phi}) to get
		\begin{align*}
			\phi (s_0 + 2^{-0\delta_0}) \leq 2^{0\delta_0} C' \phi(s_0)^{1+\delta_0} \leq 2^{0\delta_0} C' (2C')^{-1 - 1/\delta_0} \leq C'^{-1/\delta_0} 2^{-0-1-1/\delta_0}
		\end{align*}
		by the choice of $s_0$.\\
		Now, assume the claim holds true for $m \geq 0$. Then for $m+1$ using (\ref{fctl_ineq_phi}) again we get
		\begin{align*}
			\phi\left(s_0 + \sum_{i=0}^{m+1} 2^{-i\delta_0} \right) &\leq 2^{(m+1)\delta_0} C' \phi\left(s_0 + \sum_{i=0}^m 2^{-i\delta_0} \right)^{1+\delta_0}\\
			&\leq 2^{(m+1)\delta_0} C' \left(C'^{-1/\delta_0} 2^{-m-1-1/\delta_0}\right)^{1+\delta_0}\\
			&= 2^{(m+1)\delta_0} C'^{1 -1/\delta_0 -1} 2^{-m-1-1/\delta_0 -m\delta_0 -1\delta_0 -1}\\
			&= C'^{ -1/\delta_0 } 2^{-(m+1)-1-1/\delta_0 },
		\end{align*}
		where in the second inequality we apply the induction hypothesis.
	\end{proof}

	\subsection{Auxiliary Monge-Ampère equation and application} \label{subsec_auxMA_meanvalue}
	To prove the functional inequality (\ref{fctl_ineq_phi}), we will setup an auxiliary Monge-Ampère equation and show that a certain function $\Phi$ depending on its solution is non-positive. This allows to estimate quantities related to (\ref{fctl_ineq_phi}) against its solution, which is $\hat{\omega}$-psh such that an $\alpha$-invariant estimate leads to an uniform upper bound. 
	
	We fix $s\geq 0$ and assume that $\Omega_s \neq \emptyset$ since $\phi(s) = 0$ otherwise. For the setup of the Monge-Ampère equation we need to cut-off the function $v-s$ at $0$ and smooth it to obtain a smooth solution by Yau's theorem. This is done in the following Lemma.
	
	\begin{Lemma}
		There is a sequence of smooth, positive functions $\eta_k(v-s) \colon \XX \longrightarrow \RR$ converging uniformly to $\max \{0, v-s\}$ from above for $k \rightarrow \infty$.\\
		More generally, let $f \colon  \XX \longrightarrow \RR$ be continuous. Then there is a sequence of orbifold-smooth functions $\eta_k(f)$ converging uniformly to $f$. 
	\end{Lemma}
	
	\begin{proof}
		Note first that $\max \{0, v-s\}$ is continuous on $\XX$. By the construction for the general case, the non-negativity of the $\eta_k(v-s)$ will follow. We have $\|\max \{0, v-s\} - \eta_k(v-s) \|_\infty \leq \varepsilon_k$ where $\varepsilon_k$ converges to zero. Then we replace $\eta_k(v-s)$ by $\eta_k(v-s) + \varepsilon_k + 1/k$ to get strict positivity and convergence from above.\\
		In the general case, we fix a finite cover by orbifold charts $(U_\alpha, G_\alpha, \pi_\alpha)$. We fix a corresponding orbifold partition of unity $\{\varphi_\alpha \}$ that we can construct as follows: On each $U_\alpha$, choose a smooth function $\mu_\alpha$ with compact support $K_\alpha$ that is strict positive on an open set $U'_\alpha$ such that the sets $\pi_\alpha (U'_\alpha)$ still cover $\XX$. To yield an orbifold object, we need to make the $\mu_\alpha$ $G_\alpha$-invariant, for example by averaging 
		\begin{align*}
			\mu^{G}_\alpha \coloneqq  \sum_{g \in G_\alpha} \mu_\alpha \circ g.
		\end{align*}
		Note that $\mu^G_\alpha$ is still compactly supported in $U_\alpha$. Now, define
		\begin{align*}
			\varphi_\alpha \coloneqq \frac{\mu^{G}_\alpha}{\sum_\beta \mu^{G}_\beta}
		\end{align*}
		which has the desired properties of a partition of unity.\\
		Consider the convolution of the representative $f_\alpha$ of $f$ on $U_\alpha$ (we may assume that $f_\alpha$ is trivially continued on $\CC^n$) with a smooth kernel $\rho_k$. Fix a smooth, non-negative function $\rho$ with compact support in $B_1(0) \subset \CC^n$ and $L^1$-norm equal to $1$. To make the convolution $G_\alpha$ invariant we define for $x\in U_\alpha$ and $k \in \NN$
		\begin{align*}
			f_\alpha * \rho_k (x) \coloneqq \frac{1}{|G_\alpha|} \sum_{g \in G_\alpha} \int_{\CC^n} f_\alpha(y) k^n\rho(k (g(x) -y )) dy.
		\end{align*}
		Because $f_\alpha$ is continuous on $U_\alpha$ the convergence on any compact set in $U_\alpha$ is uniform. In particular, the function
		\begin{align*}
			\varphi_\alpha (f_\alpha * \rho_k )
		\end{align*}
		converges uniformly to $\varphi_\alpha f_\alpha$ on the support $K_\alpha$ of $\varphi_\alpha$ and is $G_\alpha$ invariant by construction. Therefore, it descends to $\XX$ and the sum of the local approximations will yield a global approximation as demanded.
	\end{proof}
	
	For any $r \geq 0$ we find that
	\begin{align*}
		r\phi(s+r) = (r + s -s ) \int_{\Omega_{s+r}} e^{F_t} \omega_\XX^n \leq \int_{\Omega_{s}} (v-s) e^{F_t} \omega_\XX^n \eqqcolon A_s
	\end{align*}
	so that it is enough to estimate $A_s$ to obtain (\ref{fctl_ineq_phi}). We remark that 
	\begin{align*}
		A_s \leq \frac{1}{V_t} \int_{\Omega_{s}} |v|  \omega_t^n \leq \frac{\|v\|_{L^1(\XX , \omega_t)}}{V_t} \leq 1
	\end{align*}
	by the assumption on $v$. It follows by dominated convergence that we can approximate $A_s$ by 
	\begin{align*}
		A_{s,k} \coloneqq \int_{\Omega_{s}} \eta_k(v-s) e^{F_t} \omega_\XX^n,
	\end{align*}
	in particular, $A_{s,k}\leq 2$ for $k$ large enough.
	
	We consider the auxiliary Monge-Ampère equation
	\begin{align} \label{aux_MAeq_mvineq}
		(\hat{\omega}_t + \ii \ddbar \psi_k )^n = V_t\frac{\eta_k(v-s)}{A_{s,k} }e^{F_t}\omega_\XX^n, \quad \quad \sup_\XX \psi_k = 0.
	\end{align}
	It admits an orbifold smooth solution by \cite{Faulk2019_thesis,Faulk2019_YauOrbifolds} because $\hat{\omega_t} + \ii \ddbar \varphi_t = \omega_t$ is a Kähler form and the RHS is smooth, positive and its integral over $\XX$ is equal to $V_t$.
	
	In the following, we aim to obtain an estimate for (powers of) $v-s$ and $A_{s,k}$ in terms of the solution $\psi_k$ and constants. Since $\psi_k$ is $\omega_\XX$-psh, we can then apply an $\alpha$-invariant estimate. For this, we fix two constants. First, let $C_\infty$ be the uniform $L^\infty$-bound for $\VV_t - \varphi_t$ from Section \ref{sec_Linfty_est} and define $\Lambda \coloneqq C_\infty + 1$. Next, we claim that there is $\varepsilon > 0$ such that 
	\begin{align} \label{def_vareps_mvineq}
		\varepsilon^{n+1} = \left( \frac{n+1}{n^2}\right)^n (a + n\varepsilon)^n A_{s,k}.  
	\end{align}
	To see this, define a polynomial $f \colon \RR \longrightarrow \RR$ by
	\begin{align*}
		f(x) \coloneqq x^{n+1} - \left( \frac{n+1}{n^2}\right)^n (a + nx)^n A_{s,k} = x^{n+1} - \sum_{i=0}^n A_{s,k} \left( \frac{n+1}{n^2}\right)^n \binom{n}{i} a^{n-i} n^i x^i.
	\end{align*}
	Then it is clear that $f(0) < 0$ and $\lim_{x\rightarrow \infty} f(x) = + \infty$, so by the intermediate value theorem $f$ has a positive root $\varepsilon$. Moreover, the Lagrange-Zassenhaus bound for polynomial roots yields
	\begin{align*}
		\varepsilon &\leq 2 \max_{i \in \{0, \dots , n\} } \left( A_{s,k} \left( \frac{n+1}{n^2}\right)^n \binom{n}{i} a^{n-i} n^i \right)^{\frac{1}{n-i+1}}\\
		&\leq 2 \max_{i \in \{0, \dots , n\} } \left( 2^{(\frac{1}{n-i+1}-\frac{1}{n+1})} \left( \frac{n+1}{n^2}\right)^n \binom{n}{i} a^{n-i} n^i \right)^{\frac{1}{n-i+1}} A_{s,k}^\frac{1}{n+1}
	\end{align*}
	since $A_{s,k}\leq 2$ so that we have 
	\begin{align} \label{estimate_varpeps}
		\varepsilon &\leq C(n,a) A_{s,k}^\frac{1}{n+1}
	\end{align}
	for a constant depending only on $a$ and $n$.
	
	\textbf{Claim.} The function $\Phi \colon \XX \longrightarrow \RR$, 
	\begin{align*}
		\Phi \coloneqq - \varepsilon (-\psi_{t,k} + \varphi_t + \Lambda )^\frac{n}{n+1} + v - s 
	\end{align*}
	is non-positive.
	
	\begin{Bemerkung}
		The function $\Phi$ is chosen in such a way that when computing its Laplacian, $\psi_k$ as well as $\Delta \psi_k$ appear. The latter can be estimated using the auxiliary Monge-Amp\`ere equation (\ref{aux_MAeq_mvineq}) allowing to compare the RHS of (\ref{aux_MAeq_mvineq}) and its $\omega_\XX$-psh solution $\psi_k$ directly.
	\end{Bemerkung}
	
	\begin{proof}
		First, we note that
		\begin{align*}
			-\psi_{t,k} + \varphi_t + \Lambda = (\VV_t - \psi_{t,k}) + (\varphi_t - \VV_t + C_\infty) + 1 \geq 1
		\end{align*}
		since both parenthesis are non-negative by the $L^\infty$-estimate for $\varphi_t - \VV_t$ and as $\psi_{t,k}$ is $\hat{\omega}_t$-psh.\\
		On $\XX\setminus \Omega_s$ we have $v \leq s$ and the first term of $\Phi$ is negative, so $\Phi < 0$ there. On the other hand, $\Phi$ is twice continuously differentiable in the orbifold sense on $\overline{\Omega_s}$ and for the Claim to fail it has to achieve its maximum therefore at a point $x_0 \in \Omega_s$.
		
		We apply the maximum principle in the following way. We choose a locally uniformizing system $(U, G, \pi)$ of $\XX$ around $x_0$. Then the local lift of $\Phi$ achieves its maximum at the preimages of $x_0$ and is twice continuously differentiable in a neighborhood of $\pi^{-1}(x_0)$. Because $\omega_t$ is an orbifold Kähler metric, its associated Laplace operator $\Delta_{\omega_t}$ is elliptic. The maximum principle implies
		\begin{align*}
			\Delta_{\omega_t} \Phi \leq 0
		\end{align*}
		and we can compute the Laplacian by the same local formula as in the smooth case yielding
		\begin{align*}
			\Delta_{\omega_t} \Phi &= - \varepsilon g^{j\overline{k}} \left( \frac{-n}{(n+1)^2} (-\psi_{t,k} + \varphi_t + \Lambda )^\frac{-n-2}{n+1} \left(- \frac{\partial \psi_{t,k}}{\partial z_j} + \frac{\partial \varphi_t}{\partial z_j}\right) \left(- \frac{\partial \psi_{t,k}}{\partial \overline{z}_k} + \frac{\partial \varphi_t}{\partial \overline{z}_k} \right) \right.\\
			& \quad \quad \quad \left. \quad \quad   + \frac{n}{n+1} (-\psi_{t,k} + \varphi_t + \Lambda )^\frac{-1}{n+1} \left(- \frac{\partial^2 \psi_{t,k}}{\partial z_j\partial \overline{z}_k} + \frac{\partial^2 \varphi_t}{\partial z_j\partial \overline{z}_k}\right) \right) + \Delta_{\omega_t} v.
		\end{align*}
		Because the inverse Matrix $(g^{j\overline{k}})$ of the coefficients of $\omega_t$ is positive definite, the terms in the first line of this equation give a semi-positive contribution and we can conclude using the assumption $\Delta_{\omega_t} v \geq -a$ that
		\begin{align*}
			0 &\geq \frac{\varepsilon n}{n+1} (-\psi_{t,k} + \varphi_t + \Lambda )^\frac{-1}{n+1} g^{j\overline{k}} \left(\hat{g}_{j\overline{k}} + \frac{\partial^2 \psi_{t,k}}{\partial z_j\partial \overline{z}_k} - \hat{g}_{j\overline{k}} - \frac{\partial^2 \varphi_t}{\partial z_j\partial \overline{z}_k}\right) -a\\
			&=  \frac{\varepsilon n}{n+1} (-\psi_{t,k} + \varphi_t + \Lambda )^\frac{-1}{n+1} \left(\tr_{\omega_t} (\hat{\omega}_t + \ii \ddbar \psi_{t,k}) - \tr_{\omega_t} (\hat{\omega}_t + \ii \ddbar \varphi_t ) \right) -a.
		\end{align*}
		The term $\tr_{\omega_t} (\hat{\omega}_t + \ii \ddbar \psi_{t,k})$ equals the sum of the eigenvalues $\lambda_1, \dots , \lambda_n$ of the form $\hat{\omega}_t + \ii \ddbar \psi_{t,k}$ with respect to the metric $\omega_t$. These are all positive by the Monge-Ampère equation (\ref{aux_MAeq_mvineq}) that defines $\psi_{t,k}$. Thus, we can apply the inequality 
		\begin{align*}
			\frac{\lambda_1 + \dotso + \lambda_n}{n} \geq (\lambda_1 \cdot \dotso \cdot \lambda_n)^{\frac{1}{n}}
		\end{align*}
		between the arithmetic and the geometric mean. The product of the eigenvalues of $\hat{\omega}_t + \ii \ddbar \psi_{t,k}$ on the RHS is the ratio 
		\begin{align*}
			\frac{(\hat{\omega}_t + \ii \ddbar \psi_{t,k})^{n}}{\omega_t^n} = \frac{\eta_k(v-s)}{A_{s,k}}
		\end{align*}
		by using (\ref{aux_MAeq_mvineq}) again. This gives
		\begin{align*}
			0 \geq \frac{\varepsilon n^2}{n+1} (-\psi_{t,k} + \varphi_t + \Lambda )^\frac{-1}{n+1} \left( \frac{\eta_k(v-s)}{A_{s,k}} \right)^{\frac{1}{n}} - \frac{\varepsilon n^2}{n+1} - a
		\end{align*}
		because the trace of $\omega_t$ is $n$ and $(-\psi_{t,k} + \varphi_t + \Lambda )^\frac{-1}{n+1} \leq 1$. Using $\eta_k(v-s) \geq v-s $ this simplifies further and can be rearranged to
		\begin{align*}
			a + \varepsilon n \geq \frac{\varepsilon n^2}{n+1} (-\psi_{t,k} + \varphi_t + \Lambda )^\frac{-1}{n+1} \left( \frac{v-s}{A_{s,k}} \right)^{\frac{1}{n}}
		\end{align*}
		which holds at the preimages of $x_0$ in $U$ and thus for the global objects at $x_0$. We take the $n$-th power keeping in mind that both sides are semi-positive and divide by positive constants to get
		\begin{align*}
			\left( \frac{n+1}{n^2}\right)^n (a + \varepsilon n)^n A_{s,k} (-\psi_{t,k} + \varphi_t + \Lambda )^\frac{n}{n+1} \geq \varepsilon^n (v-s).
		\end{align*}
		By the definition of $\varepsilon$ in (\ref{def_vareps_mvineq}) this is equivalent to
		\begin{align*}
			\varepsilon (-\psi_{t,k} + \varphi_t + \Lambda )^\frac{n}{n+1} \geq v -s
		\end{align*}
		at $x_0$. In other words, $\max_\XX \Phi \leq 0$, the claim is established.
	\end{proof}
	
	\subsection{Uniform \texorpdfstring{$\alpha$}{alpha}-invariant estimate}    \label{subsec_alphainv_application}
	The above claim together with the estimate (\ref{estimate_varpeps}) for $\varepsilon$ gives
	\begin{align*}
		v-s \leq C(n,a) A_{s,k}^{\frac{1}{n+1}} (-\psi_{t,k} + \varphi_t + \Lambda )^\frac{n}{n+1}
	\end{align*}
	and on $\Omega_s$ where the LHS is non-negative we can take the $\frac{n+1}{n}$-th power of this inequality so that
	\begin{align} \label{ineq_of_integrands2}
		(v-s)^\frac{n+1}{n} A_{s,k}^{-\frac{1}{n}} \leq C_1 (-\psi_{t,k} + \varphi_t + \Lambda ) \leq C_1 (-\psi_{t,k} + \Lambda ) 
	\end{align}
	where we use the normalization $\sup \varphi_t = 0$ in the second inequality.
	
	To apply an $\alpha$-invariant estimate uniformly in $t$ we choose a constant $C_2 >0 $ as detailed in \cite[p.\ 8f]{GPS22} such that $\chi \leq C_2 \omega_\XX$ and $\alpha_0 = \alpha_0(n, \omega_\XX, \chi) > 0$ such that 
	\begin{align*}
		C_1 \alpha_0 < \alpha(\XX, (C_2+1) \omega_\XX )
	\end{align*}
	where the RHS is the orbifold version of Tian's $\alpha$-invariant (cf.\ \cite{Tian87_KEMetrics} p.\ 228 and Section \ref{sec_alpha_est}). Because all $\psi_{t,k}$ are in $\Psh_{\leq 0} (\XX, (C_2+1) \omega_\XX )$, are smooth, and have supremum $0$, using (\ref{ineq_of_integrands2}), there is $C_3 >0 $ independent of $t$ and $k$ such that
	\begingroup
	\addtolength{\jot}{0.2em}
	\begin{align*} 
		\int_{\Omega_s} \exp  \left(\alpha_0 (v-s)^\frac{n+1}{n}A_{s,k}^{-1/n} \right) \omega_\XX^n  
		&\leq \int_{\Omega_s} \exp \left(\alpha_0 C_1 (-\psi_{t,k}  + \Lambda ) \right) \omega_\XX^n\\
		&= \exp(\alpha_0 C_1 \Lambda ) \int_{X} \exp \left(- \alpha_0 C_1 \psi_{t,k}  \right) \omega_\XX^n\\
		&\leq C_3.
	\end{align*}
	\endgroup
	Next, we want to get an estimate for quantities of the form $C_{s,k} (v-s)^{\frac{p(n+1)}{n}} e^{F_t}$ in terms of $\exp  \left(\alpha_0 (v-s)^\frac{n+1}{n}A_{s,k}^{-1/n} \right)$. We recall the following inequality from \cite{GPT23} which we apply at each point of $\Omega_s$:
	
	\textbf{Fact.} Let $f \in \RR_+$ be non-negative, $F \in \RR$ a real number. Then we have
	\begin{align} \label{obviousinequality}
		f^p e^F \leq e^F (1 + |F|)^p + C(p) e^{2f}
	\end{align}
	for a constant $C(p) > 0$ only depending on $p$.
	
	\begin{proof}
		There are two cases. First, if $f \leq F$ then $f^p \leq (1 + |F|)^p$ so the inequality is clearly fulfilled. In the other case, we have $e^F \leq e^f$ and $f^p \leq C(p) e^{f}$ if the constant is chosen large enough. Then the LHS is less than the second summand on the RHS.
	\end{proof}

	On $\Omega_s$ we define a function $f \coloneqq \frac{1}{2}\alpha_0 (v-s)^\frac{n+1}{n}A_{s,k}^{-1/n}$ and apply the Fact with $e^{F_t}$ and $f^p$. We have
	\begin{align*}
		\left( \frac{1}{2}\alpha_0 (v-s)^\frac{n+1}{n}A_{s,k}^{-1/n} \right)^p e^{F_t} = f^p e^{F_t} \leq e^{F_t} (1+|F_t|)^p + \exp  \left(\alpha_0 (v-s)^\frac{n+1}{n}A_{s,k}^{-1/n} \right).
	\end{align*}
	Integrating this inequality over $\Omega_s$ yields
	\begin{align*}
		\int_{\Omega_s} \left( \frac{1}{2}\alpha_0 (v-s)^\frac{n+1}{n}A_{s,k}^{-1/n} \right)^p e^{F_t} \omega_\XX^n 
		&\leq \int_{\Omega_s} e^{F_t} (1+|F_t|)^p \omega_\XX^n\\
		& \quad \quad + C(p) \int_{\Omega_s} \exp  \left(\alpha_0 (v-s)^\frac{n+1}{n}A_{s,k}^{-1/n} \right) \omega_\XX^n\\
		&\leq C(\XX, \omega_\XX, p, N) + C(p) C_3
	\end{align*}
	because the $p$-entropy of $\omega_t$ is uniformly bounded by $N$. By dividing this inequality by positive constants we get
	\begin{align*}
		\int_{\Omega_s} (v-s)^\frac{p(n+1)}{n}e^{F_t} \omega_\XX^n \leq C_4 A_{s,k}^{p/n} 
	\end{align*}
	where the RHS converges to $C_4 A_s^{p/n} $ for $k \rightarrow \infty$ but the LHS is independent of $k$.
	
	As remarked above, we want to estimate $A_s$ which we can do now by applying Hölder's inequality with exponent $\frac{p(n+1)}{n}$ and its dual exponent $p'$ to get
	\begin{align*}
		A_s = \int_{\Omega_s} (v-s) e^{F_t} \omega_\XX^n &\leq \left(\int_{\Omega_s} (v-s)^\frac{p(n+1)}{n}e^{F_t} \omega_\XX^n \right)^{\frac{p(n+1)}{n}} \left( \int_{\Omega_s} e^{F_t} \omega_\XX^n \right)^{\frac{1}{p'}}\\ 
		&\leq \left(C_4 A_s^{p/n} \right)^{\frac{p(n+1)}{n}}  \left( \int_{\Omega_s} e^{F_t} \omega_\XX^n \right)^{\frac{1}{p'}} \\
		&= C_5 A_s^{\frac{1}{n+1}} \phi(s)^{\frac{1}{p'}}
	\end{align*}
	and since we did assume $A_s > 0$ this leads to
	\begin{align*}
		r\phi(s+r) \leq A_s \leq C' \phi(s)^{\frac{n+1}{np'}}
	\end{align*}
	which is equivalent to (\ref{fctl_ineq_phi}) since $\frac{n+1}{np'} = \frac{pn +p -n}{np} = 1 + \frac{p-n}{np} \eqqcolon 1 + \delta_0$. This finishes the proof of Theorem \ref{meanvalueineq}.

	
	
	
	\section{\texorpdfstring{$L^\infty$}{L-infinity} estimate for the psh-envelope} \label{sec_Linfty_est}
	The goal of this section is to prove the $L^\infty$-estimate for $\VV_t - \varphi_t$, namely Theorem \ref{Linfty_estimate}. This will be done very similarly to the smooth case in \cite{GPTW24}. We will point out the modifications needed for the singular setting. Also, this section is somewhat parallel to the previous one. 
	
	Throughout this section, we will denote by $\Omega_s \coloneqq \{\VV_t - \varphi_t > s \}$ the upper-level set of $\VV_t - \varphi_t$. Since $\VV_t$ is the psh-envelope associated to $\hat{\omega}_t$ and $\varphi_t$ is by construction $\hat{\omega}_t$-psh we get $\varphi_t \leq \VV_t$. Thus, it is enough to show that the function
	\begin{align*}
		\phi \colon \RR_{\geq 0} \longrightarrow \RR_{\geq 0}, \quad \quad s \longmapsto \frac{1}{V_t} \int_{\Omega_s} \omega_t^n = \int_{\Omega_s} e^{F_t} \omega_\XX^n
	\end{align*}
	vanishes for $s \geq C$, where $C= C(\XX, \omega_\XX, p, n, \chi, N) >0$ only depends on the quantities stated in Theorem \ref{Linfty_estimate}.  The main part of the proof is to establish for every $s\geq 0$ and every $r \geq 0$ the functional inequality 
	\begin{align} \label{fctl_ineq_Linfty}
		r \phi(s+r) \leq C' \phi (s)^{1+\delta_0}
	\end{align}
	where $\delta_0 \coloneqq (p-n)/(np) > 0$ and $C' >0$ depends on the same quantities as $C$. However, we will allow dependence on
	\begin{align*}
		E_t \coloneqq \int_\XX (\VV_t - \varphi_t)  e^{F_t} \omega_X^n
	\end{align*}
	until we provide a bound of $E_t$ in terms of the entropy $\Ent_p(\omega_t)$ in the last part \ref{subsec_energybound} of this section. 
	
	First, we note that a DeGiorgi-type iteration yields the Theorem if (\ref{fctl_ineq_Linfty}) holds (cf.\ Section \ref{subsec_DeGiorgi} for more details):\\
	We remark that $\phi$ is not-increasing and decays at least at rate of $1/s$. Indeed, on $\Omega_s$ for $s>0$ we have $(\VV_t - \varphi_t)/s >1$ so that
	\begin{align*}
		\phi(s) < \int_{\Omega_s} \frac{\VV_t - \varphi_t}{s} e^{F_t} \omega_X^n \leq \frac{1}{s} \int_{\Omega_0} (\VV_t - \varphi_t) e^{F_t} \omega_X^n = \frac{E_t}{s}.
	\end{align*}
	We choose $s_0 \coloneqq E_t (2C')^{1/\delta_0}$ such that 
	\begin{align*}
		\phi (s_0)^{\delta_0} < \frac{1}{2C'}.
	\end{align*}
	Moreover, for $i \in \{0,1,2,\dots \}$ we define $r_i \coloneqq 2^{-i\delta_0}$ and we iterate (\ref{fctl_ineq_Linfty}) starting from $s_0$ with the $r_i$'s. We set
	\begin{align*}
		C \coloneqq s_0 + \sum_{i = 0}^{\infty} r_i = E_t(2C')^{1/\delta_0} + \frac{1}{1-2^{-\delta_0}}
	\end{align*}
	which only depends on the quantities stated in Theorem \ref{Linfty_estimate} and $E_t$. Then $\phi(s) = 0$ if $s \geq C$.

	\subsection{Auxiliary Monge-Ampère equation and application}
	As in the smooth case (cf.\ \cite{GPTW24}), the proof is based on the use of auxiliary Monge-Ampère equations, here in two instances. One is very similar to the proof of the mean-value-inequality (cf.\ Section \ref{sec_meanvalue_ineq}) which produces a $\hat{\omega}_t$-psh function such that an uniform $\alpha$-invariant estimate can be applied. However, to yield a smooth solution by the orbifold version of Yau's theorem, we need smooth coefficients. For this, we approximate the psh-envelope $\VV_t$ uniformly by the smooth functions $u_{t,\beta}$ obtained in Proposition \ref{approx_psh_envelope}.

	To assure the positivity of the RHS of the auxiliary Monge-Ampère equation we will cut off $u_{t,\beta} - \varphi_t - s$ smoothly at zero. For this, fix a sequence of positive functions $\tau_k \colon \RR \longrightarrow \RR$ such that $\tau_k (x) \geq x \chi_{\RR_+} (x) + 1/k$ and converges uniformly to $x \chi_{\RR_+} (x)$. This implies in particular that the normalizing constant
	\begin{align*}
		A_{s,k,\beta} \coloneqq \int_\XX \tau_k(-\varphi_t + u_{t,\beta} -s ) e^{F_t} \omega_\XX^n
	\end{align*}
	converges for $k, \beta \rightarrow \infty$ to 
	\begin{align*}
		A_{s} \coloneqq \int_\XX (-\varphi_t + \VV_t -s)  e^{F_t} \omega_\XX^n.
	\end{align*}
	We consider the Monge-Ampère equation
	\begin{align} \label{MAeq_Linfty_est}
		(\hat{\omega}_t + \ii \ddbar \psi_{t,k, \beta})^n = V_t \frac{\tau_k(-\varphi_t + u_{t,\beta} -s )}{A_{s,k,\beta}} e^{F_t} \omega_\XX^n, \quad \quad \sup_\XX \psi_{t,k, \beta} = 0
	\end{align}
	which admits an orbifold smooth solution because the RHS is positive and orbifold smooth and the class of $\hat{\omega}_t$ is Kähler. Because $\psi_{t,k, \beta}$ is $\hat{\omega}_t$-psh and non-positive, we have $\psi_{t,k, \beta} \leq \VV_t$ and thus, by taking $\beta$ large enough, we may assume $\psi_{t,k, \beta} \leq u_{t, \beta} +1$. 
	
	We fix the constants
	\begin{align*}
		\varepsilon \coloneqq A_{s,k,\beta}^{\frac{1}{n+1}}\left(\frac{n+1}{n}\right)^{\frac{n}{n+1}}, \quad \quad \Lambda \coloneqq \frac{n}{n+1} A_{s,k,\beta}.
	\end{align*} 
	and set up an auxiliary function
	\begin{align*}
		\Phi \coloneqq -\varepsilon (-\psi_{t,k,\beta} + u_{t,\beta} +1 + \Lambda )^{\frac{n}{n+1}} - (\varphi_t - u_{t,\beta} +s).
	\end{align*}
	
	\textbf{Claim.} There are constants $\varepsilon_\beta \geq 0$ converging to $0$ for $\beta \rightarrow \infty$ such that
	\begin{align*}
		\sup_\XX \Phi \leq \varepsilon_\beta.
	\end{align*}
	
	\begin{proof}
		First, let us note that $-\psi_{t,k,\beta} + u_{t,\beta} +1 + \Lambda > 0$ such that $\Phi$ is orbifold smooth, in particular continuous. Therefore, it attains its maximum at a point $x_0 \in \XX$.
		
		Assume first that $x_0 \in \XX \setminus \Omega_s$. Then $\VV_t > \varphi_t + s$ and we can estimate
		\begin{align*}
			\Phi (x_0) \leq - (\varphi_t - u_{t,\beta} + s) < - \VV_t + u_{t,\beta} \leq \varepsilon_\beta
		\end{align*}
		by the uniform convergence of $u_{t,\beta}$ towards the psh-envelope.
		
		On the other hand, if $x_0 \in \Omega_s$, we can apply the maximum principle $\Delta_t \Phi (x_0) \leq 0$ with the Laplacian $\Delta_t$ associated to $\omega_t$ in a local uniformizing system around $x_0$. This yields in local coordinates
		\begin{align*}
			0 &\geq - \varepsilon g^{j\overline{k}} \left(
			\frac{n}{n+1} (-\psi_{t,k,\beta} + u_{t,\beta} + \Lambda +1 )^\frac{-1}{n+1} \left(- \frac{\partial^2 \psi_{t,k,\beta}}{\partial z_j\partial \overline{z}_k} + \frac{\partial^2 u_{t,\beta}}{\partial z_j\partial \overline{z}_k}\right) \right) \\
			& \quad \;    - g^{j\overline{k}} \left(  \frac{\partial^2 \varphi_t}{\partial z_j\partial \overline{z}_k} - \frac{\partial^2 u_{t,\beta}}{\partial z_j\partial \overline{z}_k} \right)
		\end{align*}
		since we can neglect the term with the first order derivatives as in Section \ref{subsec_auxMA_meanvalue} because the inverse $(g^{j\overline{k}})$ of the coefficient matrix of $\omega_t$ is positive definite. Moreover, we add and subtract the coefficients $\hat{g}_{j\overline{k}}$ of the form $\hat{\omega}_t$ to get
		\begin{align*}
			0 &\geq  \varepsilon g^{j\overline{k}} \left( \frac{n}{n+1} ( - \psi_{t,k,\beta} + u_{t,\beta} + \Lambda +1 )^\frac{-1}{n+1} \left( \hat{g}_{j\overline{k}} + \frac{\partial^2 \psi_{t,k,\beta}}{\partial z_j\partial \overline{z}_k} - \hat{g}_{j\overline{k}} - \frac{\partial^2 u_{t,\beta}}{\partial z_j\partial \overline{z}_k}\right) \right) \\
			& \quad \;    - g^{j\overline{k}} \left( \hat{g}_{j\overline{k}} + \frac{\partial^2 \varphi_t}{\partial z_j\partial \overline{z}_k} - \hat{g}_{j\overline{k}}- \frac{\partial^2 u_{t,\beta}}{\partial z_j\partial \overline{z}_k} \right)
		\end{align*}
		which is equivalent to
		\begin{align*}
			0 &\geq \frac{\varepsilon n}{n+1} ( - \psi_{t,k,\beta} + u_{t,\beta} + \Lambda +1 )^\frac{-1}{n+1} \left(\tr_{\omega_t} (\hat{\omega}_t + \ii \ddbar \psi_{t,k,\beta}) - \tr_{\omega_t} (\hat{\omega}_t + \ii \ddbar u_{t,\beta}) \right)\\
			& \quad \; - \tr_{\omega_t} (\hat{\omega}_t + \ii \ddbar \varphi_t) + \tr_{\omega_t} (\hat{\omega}_t + \ii \ddbar u_{t,\beta}).
		\end{align*}
		As in the previous section, the term $\tr_{\omega_t} (\hat{\omega}_t + \ii \ddbar \psi_{t,k,\beta})$ equals the sum of the eigenvalues $\lambda_1, \dots , \lambda_n$ of the form $\hat{\omega}_t + \ii \ddbar \psi_{t,k,\beta}$ with respect to the metric $\omega_t$. We use the inequality 
		between the arithmetic and the geometric mean as before. Recall that the product of the eigenvalues of $\hat{\omega}_t + \ii \ddbar \psi_{t,k,\beta}$ is the ratio 
		\begin{align*}
			\frac{(\hat{\omega}_t + \ii \ddbar \psi_{t,k,\beta})^{n}}{\omega_t^n} = \frac{\tau_k(-\varphi_t +u_{t,\beta} - s)}{A_{s,k,\beta}}
		\end{align*}
		by using (\ref{MAeq_Linfty_est}). Since $\beta$ is large enough, we have  $ - \psi_{t,k,\beta} + u_{t,\beta} + \Lambda +1 > \Lambda$. Moreover, by construction, $u_{t,\beta}$ is $\hat{\omega}_t$-psh so that we have
		\begin{align*}
			0 &\geq \frac{\varepsilon n^2}{n+1} ( - \psi_{t,k,\beta} + u_{t,\beta} + \Lambda +1 )^\frac{-1}{n+1} \left(\frac{\tau_k(-\varphi_t +u_{t,\beta} - s)}{A_{s,k,\beta}}\right)^{\frac{1}{n}} - n \\
			& \quad \; + \left( 1- \frac{\varepsilon n^2}{n+1} \Lambda^{-\frac{1}{n+1} }\right) \tr_{\omega_t} (\hat{\omega}_t + \ii \ddbar u_{t,\beta})\\
			&= \frac{\varepsilon n^2}{n+1} ( - \psi_{t,k,\beta} + u_{t,\beta} + \Lambda +1 )^\frac{-1}{n+1} \left(\frac{\tau_k(-\varphi_t +u_{t,\beta} - s)}{A_{s,k,\beta}}\right)^{\frac{1}{n}} - n
		\end{align*}
		since by definition of $\varepsilon$ and $\Lambda$ the coefficient of $\tr_{\omega_t} (\hat{\omega}_t + \ii \ddbar u_{t,\beta})$ is equal to zero. We add $n$ to the resulting inequality, divide by $n$ and take the $n$-th power to obtain
		\begin{align*}
			1 &\geq \left(\frac{\varepsilon n}{n+1} \right)^n ( - \psi_{t,k,\beta} + u_{t,\beta} + \Lambda +1 )^\frac{-n}{n+1} \tau_k(-\varphi_t +u_{t,\beta} - s) A_{s,k,\beta}^{-1} \\
			&\geq \left(\frac{\varepsilon n}{n+1} \right)^n ( - \psi_{t,k,\beta} + u_{t,\beta} + \Lambda +1 )^\frac{-n}{n+1} (-\varphi_t +u_{t,\beta} - s) A_{s,k,\beta}^{-1}
		\end{align*}
		at $x_0$ because $\tau_k > \id_\RR$. This is equivalent to -- using the definition of $\varepsilon$ --
		\begin{align*}
			- (\varphi_t - u_{t,\beta} + s) \leq \varepsilon ( - \psi_{t,k,\beta} + u_{t,\beta} + \Lambda +1 )^\frac{n}{n+1},
		\end{align*}
		in other words, $\Phi(x_0) \leq 0$. This concludes the proof of the claim.
	\end{proof}

	\subsection{Uniform \texorpdfstring{$\alpha$}{alpha}-invariant estimate}
	In the previous Claim we have shown that 
	\begin{align*}
		\sup_\XX \Phi \leq \epsilon_\beta
	\end{align*}
	where the RHS tends to $0$ as $\beta \rightarrow \infty$. This is equivalent to
	\begin{align*}
		-\varphi_t + u_{t,\beta} -s \leq \varepsilon (-\psi_{t,k,\beta} + u_{t,\beta} +1 + \Lambda )^{\frac{n}{n+1}} + \epsilon_\beta.
	\end{align*}
	Since $\beta$ was chosen large enough that $\psi_{t,k,\beta} < u_{t,\beta} +1$, the RHS is positive. However, even on $\Omega_s$ the LHS may be (slightly) negative as $u_{t,\beta}$ could be less than $\VV_t$. In this case, we replace the LHS by zero. Moreover, we can replace $\Lambda$ by $A_{s,k,\beta}$ on the RHS giving
	\begin{align*}
		\max \{-\varphi_t + u_{t,\beta} -s , 0 \}  \leq \varepsilon (-\psi_{t,k,\beta} + u_{t,\beta} +1 + A_{s,k,\beta} )^{\frac{n}{n+1}} + \epsilon_\beta.
	\end{align*}
	Taking the $\frac{n+1}{n}$-th power leads to
	\begingroup
	\addtolength{\jot}{0.7em}
	\begin{align*}
		\max \{-\varphi_t + u_{t,\beta} -s , 0 \}^\frac{n+1}{n}  &\leq \left(\varepsilon (-\psi_{t,k,\beta} + u_{t,\beta} +1 + A_{s,k,\beta} )^{\frac{n}{n+1}} + \epsilon_\beta\right)^\frac{n+1}{n}\\
		&\leq 2^{1/n} \varepsilon^\frac{n+1}{n} (-\psi_{t,k,\beta} + u_{t,\beta} +1 + A_{s,k,\beta} ) + 2^{1/n} \epsilon_\beta^\frac{n+1}{n}\\
		&\leq 2^{1/n} \frac{n+1}{n} A_{s,k,\beta}^{1/n} (-\psi_{t,k,\beta}  +1 + A_{s,k,\beta} ) + \epsilon'_\beta
	\end{align*}
	\endgroup
	where the second inequality uses the estimate $(a+b)^r \leq 2^{r-1}(a^r + b^r) $ for $r \geq 1$ and $a,b \in \RR_+$ and in the last step we use the definition of $\varepsilon$ and that $u_{t,\beta} \leq 0$. We conclude altogether that
	\begin{align} \label{ineq_of_integrands}
		\frac{\max \{-\varphi_t + u_{t,\beta} -s , 0 \}^\frac{n+1}{n}}{A_{s,k,\beta}^{1/n}}  \leq C_n (-\psi_{t,k,\beta}  +1 + A_{s,k,\beta} ) + \epsilon'_\beta.
	\end{align}
	
	As in the previous section and detailed in \cite[p.\ 8f]{GPS22}, choose a constant $C_1 >0 $ such that $\chi \leq C_1 \omega_\XX$ and $\alpha_0 = \alpha_0(n, \omega_\XX, \chi) > 0$ such that 
	\begin{align*}
		C_n \alpha_0 < \alpha(\XX, (C_1+1) \omega_\XX )
	\end{align*}
	where the RHS is the orbifold version of Tian's $\alpha$-invariant (cf.\ \cite{Tian87_KEMetrics} p.\ 228 and Section \ref{sec_alpha_est}). In this way, we can apply Hörmander's and Tian's estimate invariantly of $t$ because all $\psi_{t,k,\beta}$ will be in $\Psh_{\leq 0} (X, (C_1+1) \omega_\XX )$. Using (\ref{ineq_of_integrands}), there is $C_2 >0 $ independent of $t$, $\beta$, and $k$ such that
	\begingroup
	\addtolength{\jot}{0.7em}
	\begin{align*} 
		\int_{\Omega_s} \exp & \left(\alpha_0 \frac{\max \{-\varphi_t + u_{t,\beta} -s , 0 \}^\frac{n+1}{n}}{A_{s,k,\beta}^{1/n}}\right) \omega_\XX^n  \\ &\leq \int_{\Omega_s} \exp \left(\alpha_0 C_n (-\psi_{t,k,\beta}  +1 + A_{s,k,\beta} ) + \epsilon'_\beta\right) \omega_\XX^n\\
		&= \exp(\alpha_0 (C_n (A_{s,k,\beta} +1) + \epsilon'_\beta)) \int_{X} \exp \left(- \alpha_0 C_n \psi_{t,k,\beta}  \right) \omega_\XX^n\\
		&\leq C_2 e^{C(A_{s,k,\beta} + \epsilon'_\beta)}.
	\end{align*}
	Now, we take the limit $\beta \rightarrow \infty$. To this extend, note that $A_{s,k,\beta} + \epsilon'_\beta$ converges to $A_{s,k}$ as $u_{t,\beta}$ converges to $\VV_t$ uniformly. In particular, this gives a bound for the limes inferior of the LHS. Moreover, the integrands on the LHS are non-negative functions which converge pointwise to 
	\begin{align*}
		\exp  \left(\alpha_0 \frac{(-\varphi_t + \VV_t -s )^\frac{n+1}{n}}{A_{s,k}^{1/n}}\right) \omega_\XX^n.
	\end{align*}
	Hence, Fatou's Lemma provides the estimate 
	\begin{align*} 
		\int_{\Omega_s} \exp & \left(\alpha_0 \frac{(-\varphi_t + \VV_t -s )^\frac{n+1}{n}}{A_{s,k}^{1/n}}\right) \omega_\XX^n  \leq C_2 e^{CA_{s,k} }.
	\end{align*}
	Now, we take the limit $k \rightarrow \infty$. By the uniform convergence of $\tau_k$ we have that $A_{s,k} \rightarrow A_s$ fot $k \rightarrow \infty$. Moreover, on the RHS we can estimate
	\begin{align*}
		A_s = \int_{\Omega_s} (\VV_t -\varphi_t - s) e^{F_t} \omega_\XX^n \leq \int_\XX (\VV_t - \varphi_t) e^{F_t} \omega_\XX^n = E_t
	\end{align*}
	giving that 
	\begin{align*}
		\int_{\Omega_s} \exp & \left(\alpha_0 \frac{(-\varphi_t + \VV_t -s )^\frac{n+1}{n}}{A_{s}^{1/n}}\right) \omega_\XX^n  \leq C_2 e^{CE_t }.
	\end{align*}
	\endgroup

	We apply the inequality (\ref{obviousinequality}) from Section \ref{subsec_alphainv_application} pointwise on $\Omega_s$ to $F_t$ and
	\begin{align}
		f \coloneqq \frac{\alpha_0 (-\varphi_t + \VV_t -s )^\frac{n+1}{n}}{2 A_{s}^{1/n}} \geq 0
	\end{align}
	which yields by integration
	\begingroup
	\addtolength{\jot}{0.2em}
	\begin{align*}
		\int_{\Omega_s } \left( \frac{\alpha_0 (-\varphi_t + \VV_t -s )^\frac{n+1}{n}}{2 A_{s}^{1/n}} \right)^p e^{F_t} \omega_\XX^n &\leq  \int_{\Omega_s } e^{F_t} (1+ |F_t|)^p \omega_\XX^n\\
		& \quad + \int_{\Omega_s} \exp  \left(\alpha_0 \frac{(-\varphi_t + \VV_t -s )^\frac{n+1}{n}}{A_{s}^{1/n}}\right) \omega_\XX^n\\
		&\leq 2^{p-1} (\Ent_p(\omega_t) + 1) + C_2 e^{CE_t}.
	\end{align*}
	This is equivalent to 
	\begin{align*}
		\int_{\Omega_s } (-\varphi_t + \VV_t -s )^{p\frac{n+1}{n}} e^{F_t} \omega_\XX^n &\leq \left( \frac{2}{\alpha_o}\right)^p A_s^{\frac{p}{n}}\left(2^{p-1} (\Ent_p(\omega_t) + 1) + C_2 e^{CE_t} \right)\\
		&\leq C(\XX, \omega_\XX, p, n, N, \chi, E_t) A_s^{\frac{p}{n}}
	\end{align*}
	which we use to estimate $A_s$ by Hölder inequality with exponent $\frac{p(n+1)}{n}$ and its dual exponent $p'$
	\begin{align*}
		A_s &= \int_{\Omega_s }  (-\varphi_t + \VV_t -s )e^{F_t} \omega_\XX^n\\
		&\leq \left(\int_{\Omega_s } (-\varphi_t + \VV_t -s )^\frac{p(n+1)}{n} e^{F_t} \omega_\XX^n\right)^{\frac{n}{p(n+1)}} \left(\int_{\Omega_s } e^{F_t} \omega_\XX^n \right)^{\frac{1}{p'}}\\
		&\leq C_3 A_s^{\frac{1}{n+1}} \phi(s)^{\frac{1}{p'}}.
	\end{align*}
	This implies for any $r \geq 0$
	\begin{align*}
		r \phi(s+r) \leq A_s \leq C' \phi(s)^{\frac{1+n}{np'}} = C' \phi(s)^{1 + \frac{p-n}{np}}
	\end{align*}
	\endgroup
	which is the functional inequality (\ref{fctl_ineq_Linfty}) we had to prove.
	
	\subsection{Bound of \texorpdfstring{$E_t$}{E} by the entropy} \label{subsec_energybound}
	The constant $C'$ depends up to now on the term $E_t= \int_\XX (\VV_t - \varphi_t) e^{F_t} \omega_\XX^n$. We recall first that $\VV_t$ is non-positive giving
	\begin{align*}
		E_t \leq \int_\XX -\varphi_t e^{F_t} \omega_\XX^n.
	\end{align*}
	The function $\varphi_t$ is $\hat{\omega}_t$-psh, in particular $(C_1 +1)\omega_\XX$-psh. This allows to apply an $\alpha$-invariant estimate uniform in $t$. That is, there are constants $\alpha, C > 0$ independent of $t$ such that
	\begin{align*}
		\int_\XX \exp(-F_t - \alpha \varphi_t) e^{F_t} \omega_\XX^n = \int_\XX e^{ - \alpha \varphi_t} \omega_\XX^n \leq C.
	\end{align*}
	By the normalization $\int_\XX e^{F_t}\omega_\XX^n = 1$, we can apply Jensen's inequality for concave functions with the logarithm yielding
	\begin{align*}
		\int_\XX (-F_t - \alpha \varphi_t) e^{F_t} \omega_\XX^n &\leq \log \left( \int_\XX \exp(-F_t - \alpha \varphi_t) e^{F_t} \omega_\XX^n \right) \\
		&\leq \log C
	\end{align*}
	which implies
	\begin{align*}
		\int_\XX -\varphi_t e^{F_t} \omega_\XX^n &\leq \frac{1}{\alpha} \int_\XX F e^{F_t} \omega_\XX^n + \frac{\log C}{\alpha}\\
		&\leq C(\XX, \omega_\XX, \chi) \int_\XX |F| e^{F_t} \omega_\XX^n + C(\XX, \omega_\XX, \chi)\\
		&\leq C(\XX, \omega_\XX, \chi, p) \Ent_p(\omega_t) + C(\XX, \omega_\XX, \chi)
	\end{align*}
	by Hölder inequality. This makes the reasoning in this section independent of $E_t$ because $\Ent_p(\omega_t) \leq N$ by assumption.

	\bigskip

	\printbibliography[title={Bibliography}]

	\pagestyle{empty}


\end{document}